\documentclass[11pt]{article}

\usepackage[T1]{fontenc}
\usepackage[utf8]{inputenc}
\usepackage[frenchb]{babel}
\usepackage[bookmarks,hidelinks,pdftex]{hyperref}
\usepackage[usenames, dvipsnames]{xcolor}
\usepackage{graphicx}
\DeclareGraphicsExtensions{.jpg}
\usepackage{tikz}
\title{Some facts about the life and the scientific work of the Belgian mathematician Paul Mansion (1844-1919)
Centenary Paul Mansion ({\it working paper, 2019})\\
Quelques éléments sur la vie et l'oeuvre scientifique du mathématicien belge Paul Mansion (1844-1919)\\
Centenaire Paul Mansion ({\it document de travail, 2019})}
\author{Hervé Le Ferrand\footnote{Institut de Mathématiques de Bourgogne, leferran@u-bourgogne.fr}}
\date{\today}

\begin{document}

\maketitle
%\markboth{}

In this article, we are interested in the life and scientific work of the Belgian mathematician Paul Mansion. The year 2019 marks the centenary of his passing. We bring some new insights into Paul Mansion's work thanks to his scientific correspondence at the Royal Library in Brussels.
 
Nous nous intéressons dans cet article à la vie et à l'oeuvre scientifique du mathématicien belge Paul Mansion dont l'année 2019 marque le centenaire de sa disparition. Nous apportons quelques éclairages nouveaux sur ses travaux grâce notamment à sa correspondance scientifique qui se trouve à la bibliothèque royale à Bruxelles.

\newpage

%\newpage
%\tableofcontents
\newpage

\section{Introduction}
Si l'on devait dire en une phrase qui est Paul Mansion, peut-être écrirait-on que c'est un mathématicen et historien des sciences belge qui a vécu de 1844 à 1919 et qui à vingt trois ans débute une longue carrière à l'Université de Gand. Mais Paul Mansion est évidemment plus que cela. Son oeuvre et son rayonnement scientifiques sont exceptionnels. On ne peut qu'être impressionné par sa correspondance scientifique qui témoigne d'une activité inlassable. Mansion, l'homme, est remarquable, comme nous les verrons, par ses questionnements, ses centres d'intérêts et ses différentes activités.

Paul Mansion a écrit plus de trois cent articles dans plusieurs domaines comme ceux des Probabilités et Statistique, Equations aux dérivées partielles, Géométrie, ou encore celui de l'Histoire des Mathématiques. Nous allons parcourir l'oeuve et la vie de ce grand savant disparu il y a un siècle.

%Ce sont les écrits, mais aussi l'activité scientifique de façon générale, de Paul Mansion en Probabilité %et Statistique qui vont être l'objet de notre étude. Nous chercherons justement à comprendre %pourquoi Paul Mansion s'est préoccupé de Probabilités, quelles influences il a subies et en retour %quelle portée ont eu ses travaux. En fait quelle est la place de Paul Mansion \og probabiliste-%statisticien\fg\ dans l'école, que nous décrirons, belge de probabilités-statistique ?

\section{Quelques éléments biographiques}

Paul Mansion écrivait et recevait beaucoup de lettres. Ces documents ont été en partie conservés dans plusieurs fonds d'archives. Citons le fonds Paul Mansion de la Bibliothèque Royale de  Belgique à Bruxelles\footnote{En 1961, Lucien Kieffer mentionne dans un article sur le mathématicien Isidore Clasen (1829-1902) l'existence d'un recueil de correspondance passive de Paul Mansion acquis par la Bibliothèque Royale de Belgique en 1959. Au début de l'année 2016, nous avons contacté le service des Manuscrits de la Bibliothèque Royale de Belgique (BR) à Bruxelles. On nous a alors indiqué que la BR possède en effet deux fardes de correspondance adressée à Paul Mansion, plus de 300 documents. 
Un projet de numérisation était alors lancé avec pour objectif la mise en ligne sur le site de la BR des numérisations. En six mois seulement, le projet a été réalisé. Il est essentiellement question de mathématiques dans ces lettres. Il y souvent de longs développements. Les sujets sont très variés. Donnons quelques noms : Catalan, Clebsch, Darboux, Dedekind, Dumoulin, Hermite, Jordan, Killing, Kowalesky, Kronecker, De La Vallée Poussin, Sophus Lie, Lechalas, Mittag-Leffler, Neuberg, Wassilieff, Turquem. }, le fonds Pierre Duhem de l'Académie des Sciences de Paris\footnote{Lettres à Pierre Duhem.}, le fonds Paul Mansion de l'Université de Gand\footnote{Correspondance privée, lettres de Paul Mansion à Charles Lagasse, manuscrits de cours.} et le dossier Paul Mansion à l'Académie Royale de Belgique. Ajoutons à cela différentes notices et éloges sur Paul Mansion et nous avons un corpus important de documents permettant de cerner le personnage. Donnons tout d'abord quelques éléments biographiques. Nous nous appuyons ici essentiellement sur l'éloge de Paul Mansion faite par un de ses anciens élèves et collègue Alphonse Demoulin\footnote{Alphonse Demoulin étudia et enseigna à l'Université de Gand à partir de 1893 (jusqu'en 1936)\cite{doc2}. Il fut aussi l'élève de Gaston Darboux (1842-1917) à Paris en 1892. Demoulin travailla essentiellement dans le domaine de la Géométrie Différentielle et ses travaux furent récompensés par troix prix de l'Académie des Sciences de Paris : le prix du baron de Joest en 1906, le prix Bordin en 1911\footnote{Sujet proposé : \og {\it perfectionner en un point important la théorie des systèmes triples de surfaces orthogonales} \fg.} et le prix Poncelet en 1945. } (1869-1947) en 1927 devant la Classe des Sciences de l'Académie royale des sciences, des lettres et des beaux-arts de Belgique\footnote{Cette académie est divisée en quatre classes.} \cite{Demoulin} et sur l'article paru dans le bulletin du Cercle d'Histoire et de Folklore de Marchin \cite{Marchin}.

Paul Mansion est né le 3 Juin 1844 à Marchin-Lez-Huy dans la province de Liège, neuvième enfant d'une famille de dix enfants. Paul Mansion vient donc au monde dans une Belgique dont les frontières ont été fixées en 1839\footnote{La Belgique acquiert son indépendance en 1830. Le 21 Juillet 1831, Léopold Ier (1790-1885) monte sur le trône de Belgique.} par le traité des XXIV articles. Ses parents demeurent au hameau de Bel-Air, au lieu-dit Belgrade. Paul-Joseph Mansion, son père est receveur communal alors âgé de cinquante-trois ans. La mère de Paul Mansion, Marie-Elisabeth-Ferdinande née Deveux, âgée de trente-neuf ans, veille sur le foyer. La famille Mansion est relativement aisée, propriétaires de différents biens immobiliers et de terres \cite{Marchin}. Paul Mansion est orphelin de père dès son enfance puis orphelin de mère à dix-sept ans. Il pourra s'appuyer sur ses frères aînés. 

Quand Paul Mansion entre à l'école communale de Marchin,quelques années après la promulgation de loi organique dite \og loi Nothomb\fg\footnote{Jean-Baptiste Nothomb (1805-1881), libéral, chef du gouvernement belge de 1841 à 1845.} en 1842. Comme l'indiquent les auteurs de  \cite{Ligue}, cette loi \og oblige les communes à entretenir une école, cette dernière pouvant être subsidiée ou adoptée.\fg. Cette même loi stipule que les instituteurs doivent être diplomés et que les diplômes peuvent être délivrés par une école normale d'Etat ou privée agréée par l'Etat \cite{Ligue}. Paul Mansion a d'ailleurs comme maître Jean-Joseph Blaise qui entre dans cette catégorie d'instituteurs diplômés \cite{Marchin}. Paul Mansion se révèle un élève brillant. Il continue sa scolarité  à Huy.

Marchin-Lez-Huy fait partie  de la province de Liège mais c'est à Gand cependant que Paul Mansion effectue ses études supérieures. Il entre en 1862 à l'Ecole normale des sciences qui dépend de l'Université de Gand et en 1865, il en sort {\it professeur agrégé de l'enseignement moyen du degré supérieur pour les sciences} et dès l'automne 1865, il enseigne à l'Ecole préparatoire du Génie civil à Gand.  

\section{Formation}

Le 13 Août 1867, devant un jury constitué de professeurs des universités de Gand et Bruxelles, Paul Mansion devient {\it docteur en sciences physiques et mathématiques}. A cette époque, Mansion est au contact des deux mathématiciens Félix Dauge (1829-1889) et Mathias Schaar (1817-1867). Ce dernier enseigne l'arithmétique et suite à son décès, Paul Mansion est nommé sur la chaire de Calcul Différentiel et Intégral, et d'Analyse Supérieure de l'Université de Gand le 3 Octobre 1867 à l'âge de 23 ans. Quant à Dauge, c'est son cours de  {\it méthodologie de mathématique} qui marque durablement Mansion.

L'année 1868 voit la publication du permier article de Paul Mansion,  un article de probabilités intitulé \og {\it Sur le problème des partis} \fg\ publié dans les Mémoires de l'Académie royale des sciences de Belgique. Nous y reviendrons. Mais c'est plutôt dans le domaine de l'Analyse que Mansion va poursuivre ses recherches.

Revenons sur les mathématiciens qui ont participé au développement des mathématiques à l'Université de Gand. Nous avons déjà cité Félix Dauge et Mathias Schaar. Une figure importante des mathématiques gantoises est Jean-Guillaume Garnier (1766-1840). Mansion écrit dans le Liber Memorialis au sujet de Garnier :
\begin{quotation}
\begin{it}
Son enseignement
était diffus et ennuyeux ; en revanche, il avait avec ses
élèves des relations personnelles et sa conversation primesautière,
claire et démonstrative, était éminemment excitatrice
et pleine de renseignements historiques précis sur les récents
progrès des sciences. Aussi, on peut dire de lui qu'il a été
le principal rénovateur de l'étude des hautes mathématiques
en Belgique, par ses ouvrages et surtout par ses élèves. Parmi
ceux-ci nous citerons Quetelet, Timmermans, Verhulst,
Lemaire, Ed. Lannoi, L. Casterman, A. Leschevin, Mareska,
Ch. Morren, E. Manderlier, Fr. Duprez, A. Goethals. Garnier
était à peu près le seul professeur de la Faculté des sciences
de Gand qui ne fit pas ses leçons en latin. Il fut l'un des
fondateurs et collaborateurs des deux recueils savants de
l'époque, les Annales belgiques et la Correspondance mathématique et physique. Il appartenait à l'Académie royale de Bruxelles depuis le 7 mai 1818.
\end{it}
\end{quotation}
Plusieurs élèves de Garnier poursuivront une carrière en Mathématiques. Le plus célèbre est certainement Alphonse Quetelet (1796-1874)\footnote{Pour des analyses des travaux de Quetelet, on pourra se reporter notamment aux analyses de Jean-Jacques Droesbeke comme \cite{JJD1} ou \cite{JDD2}.}. Jean-Guillaume Garnier encadre la thèse Quetelet que ce dernier soutient le 24 Juillet 1819 à Gand. Les deux hommes fondent en 1825 le périodique {\it Correspondance mathématique et physique}\footnote{Voir le site de l'Institut royal météorologique belge :

 http://www.meteo.be/meteo/view/fr/6044768-Histoire.html}.

Sur Jean-François Lemaire (1797-1852), on lit dans le Liber Memorialis :
\begin{quotation}
\begin{it}
Quelques discours
insérés dans les Annales des Universités de Belgique, quelques
travaux statistiques dans la Correspondance mathématique
de M. Quetelet complètent à peu près son bilan. Lemaire
était avant tout un excellent professeur, un homme instruit
et d'un commerce agréable; peut-être se serait-il élevé plus
haut s'il lui avait été donné de suivre sa vocation.
\end{it}
\end{quotation}

Quant à Eloi-Joachim-Joseph Manderlier (1795-1884), il fut déchargé en 1863 du cours d'Algèbre, cours qu'assuma alors Félix Dauge. 

Jean-Timmermans (1801-1864) enseigna à Gand et ses recherches portaient sur les équations aux dérivées partielles et les équations différentielles. Il fut un des professeurs de Mansion. Mansion écrit dans Liber Memorialis :
\begin{quotation}
\begin{it}
Le 12 octobre 1833, il avait été
élu membre de l'Académie, qui l'avait couronné en 1831
pour ses Recherches sur la forme la plus avantageuse à
donner aux ailes des moulins à vent (Mém. cour. in-4,
1831, t. VIII, 26 pp. in-4). On peut citer parmi ses travaux
académiques un Mémoire sur les solutions singulières des
équations différentielles (Mém. in-4, t. XV, 1842; 24 pp.),
ses Recherches sur les axes principaux d'inertie et sur les
centres de percussion (Ibid. t. XXI, 1847, 33 pp.), et un Mémoire
sur l'intégration des équations linéaires aux dérivées partielles
à coefficients variables (Ibid. t. XXVIII, 1854, 10 pp.).
II a aussi inséré dans les Bulletins de l'Académie de très
nombreux rapports qui témoignent de l'étendue de ses connaissances;
puis, quelques notes originales dont l'une Sur le
parallélogramme des forces de Simon Stévin [Bulletins,
1846, 7 p. t. XIII) mérite d'être signalée, parce que l'auteur en
déduit le principe des vitesses virtuelles.
\end{it}
\end{quotation}

Pierre-François Verhults (1804-1849), un autre nom célèbre,  \og{\it cultiva d'abord les probabilités sous la direction de Quetelet et enseigna l'analyse au musée de Bruxelles, puis, à partir de 1834, à l'Ecole militaire}\fg  comme l'écrit Maurice Alliaume (1882-1931) dans \cite{Alliaume1930}. Verhults poursuit des travaux en Statistique : en 1844 (et 1846) il publie le mémoire {\it Recherches mathématiques sur la loi d'accroissement de la population} dans lequel il introduit la fameuse {\it équation logistique}.

Revenons sur le parcours de Jean-François Lemaire. Celui-ci quitte l'Université de Gand en 1830\footnote{La faculté des Sciences de Gand est fermée un temps après la révolution de 1830.} pour celle de Liège. En 1847, Antoine Meyer (1803-1857) lui succède à la chaire d'Analyse mathématique de l'Université de Liège. Dans \cite{Breny}, Henri Breny analyse la contribution de Meyer au calcul des Probabilités. L'oeuvre majeure de Meyer est son livre, publié après sa mort, {\it Essai sur une exposition nouvelle de la théorie des probabilités à postériori}. Breny souligne :
\begin{quotation}
\begin{it}
On constate donc que cette oeuvre de grande ampleur ne contient rien de vraiment neuf quant à la pensée probabiliste : son originalité réside dans la rigueur que met l'auteur dans l'étude de la convergence des séries et des intégrales qu'il y utilise et dans son habileté à simplifier de longues démonstrations analytiques.
\end{it}
\end{quotation}
Le livre est d'ailleurs traduit en allemand sous la direction d'Emanuel Czuber (1851-1925). Dans \cite{Breny}, est relaté les \og démêlés \fg\ de Meyer au sein de l'Académie des Sciences -dont il était membre-pour la publication de ses derniers travaux en Probabilités\footnote{Ce conflit a pour origine une critique par Meyer de travaux de Jean-Baptiste Liagre (1815-1891), capitaine à l'époque-il deviendra général-, collaborateur de Quetelet à l'Observatoire.}. Son adversaire en l'occurence fut Schaar (professeur à Gand).

\section{Thèse, mariage, fondation de Mathesis}

En 1870, Paul Mansion soutient une thèse intitulée {\it Théorie de la multiplication des fonctions elliptiques} à l'Université de Liège. Il passe ensuite un semestre en Allemagne à l'Université de 
G\"oettingen au contact de Alfred Clebsch (1833-1872) et Ernst Schering (1833-1897). A son retour, il se fiancie avec Cécile Belpaire. A son grand ami Charles Lagasse\footnote{Le chevalier Charles Lagasse de Locht (1845-1937), ingénieur des Ponts et Chaussées, secrétaire général au Ministère des Travaux publics.}\footnote{Lettre du 10/6/1871, fonds Mansion, Université de Gand.}, il écrit en 1871 : 
\begin{quote}
\begin{it}
Oublie toutes mes contradictions touchant Mademoiselle Belpaire. Tout est fini. Je suis fiancé à Melle Cécile depuis jeudi.
\end{it}
\end{quote}
Le père de la fiancée de Paul Mansion est Alphonse Belpaire\footnote{Son père Antoine Belpaire est un ancien élève de l'Ecole Polytechnique. Le second fils de Antoine Belpaire, donc le frère de Alphonse Belpaire, est Alfred Belpaire, ingénieur mécanicien diplômé de l'Ecole Centrale des Arts et Manufacture de Paris \cite{doc4}. Il créa et fit fabriquer en Belgique plusieurs modèles de locomotives.} (1817-1854), ingénieur. Alphonse Belpaire dirigea notamment le service de l'entrepôt et de la station commerciale d'Anvers, et le service spécial du Rupel ( voir la notice de Alphonse Belpaire dans \cite{doc4} et celle de Eugène Boudin dans \cite{doc3}. Cécile Belpaire a pour soeur l'écrivain Maria Belpaire (1853-1948). L'épouse de Alphonse Belpaire est Betsy Teichmann (1821-1900) fille du gouverneur de la Province d'Anvers, Theodoor Teichmann.
De cette union, naîtront : Joseph Mansion (1877-1937), linguiste à l'Université de Liège, Elisabeth Mansion, Augustin Mansion (1882-1966), prêtre, professeur de philosophie à l'Université Catholique de Louvain, Hubert Mansion (1883-1960), industriel.

En dehors de ses activité académiques, Paul Mansion participe à des cercles de réflexion comme le \og cercle Leibnitz\fg\ ou le cercle \og cercle Cauchy\fg créé par son ami Lagasse. Au père de celui-ci, il écrit en 1870\footnote{Fonds Mansion de l'Université de Gand.} :
\begin{quote}
\begin{it}
Je prépare une causerie sur l'esprit des bêtes pour le cercle Leibnitz de Gand et votre société d'ouvriers. Je m'en vais d'abord lire quelques philosophes là dessus pour connaitre les principes, puis je dépouillerai les naturalistes, et j'espère qu'il ressortira de ces études un entretien sérieux pour Gand, un amusant pour Nivelles.
\end{it}
\end{quote}
Signalons que Paul Mansion est proche des jésuites. Il est membre de la congrégation de Notre-Dame des Septs Douleurs, membre de la Saint Vincent de Paul gantoise. Le prêtre Désiré Mercier(1851-1926), professeur de philosophie à l'université catholique de Louvain\footnote{Désiré Mercier fonde en 1889 à l'Université Catholique de Louvain, l'Institut Supérieur de Philosophie.}, nommé cardinal en 1907,  est un ami de Mansion.

En 1873, Paul Mansion remporte un prix de l'Académie des Sciences de Belgique dont le sujet est : {\it résumer et simplifier la théorie des équations aux dérivées partielles des deux premiers ordres}. Cet important mémoire paraît sous forme d'un livre en 1875\footnote{Edité par Gauthier-Villars, Paris.} dans lequel Mansion fait souvent référence aux travaux novateurs de Sophus Lie (1842-1899).

En 1874, avec Eugène Catalan (1814-1894)\footnote{Sur la vie et l'oeuvre de Catalan voir \cite{doc5}.}, alors professeur à l'Université de Liège, Paul Mansion fonde le journal {\it Nouvelle correspondance mathématique}. Alphonse Demoulin (1869-1947)\cite{Demoulin} précise :
\begin{quotation}
\begin{it}
(...) ce journal ne devait s'occuper que
des parties de la science mathématique enseignées
dans les classes supérieures des établissements
d'instruction moyenne et dans les Ecoles d'ingénieurs.
Mais Catalan, qui le dirigea seul, de 1876
à 1880, en éleva peu à peu le niveau, de sorte
qu'il répondit de moins en moins aux besoins
de l'enseignement en Belgique.
\end{it}
\end{quotation}

Catalan correspond en fait avec Mansion dès l'année 1869. Dans une lettre\footnote{Fonds Mansion de la Bibliothèque Royale de Belgique.} 20/10/1869, Catalan écrit :
\begin{quote}
\begin{it}
En attendant, si le voyage de Gand à Liège, ne vous parait pas trop long, venez me trouver un dimanche (en m'avertissant afin que vous me rencontriez chez moi) : j'aurai grand plaisir à faire connaissance avec vous et à vous serrer la main. Donc, à bientôt ; à dimanche prochain par exemple.

Votre ancien

Catalan
\end{it}
\end{quote}
Paul Mansion écrit  à Charles Lagasse\footnote{Fonds Mansion de l'Université de Gand.}  en Septembre 1870 à propos de Catalan : 
\begin{quote}
\begin{it}
A Liège, j'ai vu M. Catalan que les misères de la France affligent profondément (...) Il parait qu'à Liège, et à Verviers, l'opinion publique n'est pas favorable aux Prussiens.
\end{it}
\end{quote}

Dès 1881, cette revue journal est remplacée par une nouvelle publication, {\it Mathesis}, créée par Mansion et Joseph Neuberg (1840-1926). Concernant {\it Mathesis}, Demoulin \cite{Demoulin} indique :
\begin{quotation}
\begin{it}
Son programme
était le même que le programme primitif
de la Nouvelle Correspondance. Mansion et Neuberg
eurent la sagesse de maintenir Mathesis au
niveau qu'ils lui avaient assigné et leur journal
put vivre et prospérer. Ils l'ont publié sans interruption
jusqu'à la fin de 1915, soit pendant trente cinq
années. Seules les difficultés nées de la
guerre les ont forcés à en suspendre la publication.
Mansion ne devait pas la reprendre. Mathesis
n'a recommencé à paraître qu'en 1922, sous la
direction de Neuberg et de M. Mineur.
\end{it}
\end{quotation}
La revue cesse de paraître en 1962.

\section{Un professeur influent}
Par la suite, de 1880 à 1886, Paul Mansion se consacre à des questions de {\it quadratures approchées}. Puis en 1892, il succède à Emmanuel-Joseph Boudin (1820-1893) à la chaire de calcul des Probabilités de l'Université de Gand. Emmanuel-Joseph Boudin est ingénieur de formation, diplômé de l'Ecole du Génie Civil de Gand. Il fait partie du corps des Ponts et Chaussées et participe notamment à la construction de plusieurs écluses. En Décembre 1846, il est chargé d'un cours de Probabilités à l'Université de Gand et à l'Ecole du Génie Civil. Paul Mansion a rédigé sa notice biographique parue dans le Liber Memorialis, Université de Gand (1913). Mansion écrit :
\begin{quotation}
\begin{it}
Quant au Cours de calcul des probabilités, tout imprégné
des idées les meilleures de Laplace, c'est un vrai chef-d'oeuvre
sous le rapport des principes et de l'ordre des matières, supérieur
aux meilleurs manuels. La théorie des erreurs y repose
sur l'hypothèse de Hagen\footnote{Gotthilf Heinrich Ludwig Hagen (1797 -1884), ingénieur allemand.} dont Boudin, le premier et longtemps
le seul, avait reconnu toute la valeur philosophique, s'écartant
avec raison de Laplace sur ce point. L'auteur de cette notice
espère quelque jour s'acquitter d'une dette de reconnaissance
envers son ancien maître en publiant une édition définitive et
un peu rajeunie au point de vue analytique de ce beau cours.
Boudin lui en avait donné l'autorisation, quelques années
avant sa mort.
\end{it}
\end{quotation}

Quelques année auparavant, Paul Mansion avait pris la charge du cours d'Histoire des sciences mathématiques à l'Ecole normale des sciences, cours qu'il continuera d'assurer après la disparition de cette école, dans le cadre du doctorat en sciences physiques et mathématiques. 

En Novembre 1910, Arthur Claeys (1875-1949) succède à Mansion dans les cours de Probabilités et d'Histoire des sciences. Cependant Claeys n'a pas fait de recherche en Probabilités. Mansion eut aussi comme élève, puis collègue, Junius Massau (1852-1909), ingénieur et géomètre qui développa des techniques d'{\it intégration graphique}.

Deux anciens élèves de Paul Mansion vont devenir des historiens des sciences célèbres, le prêtre Henri Bosmans (1852-1928) et Georges Sarton (1884-1956). Ajoutons à cela que l'enseignement d'histoire des sciences de Bosmans marquera profondément Georges Lemaître (1894-1966), un des père du Big Bang. Georges Lemaître a eu en effet Henri Bosmans comme professeur au Collège (jésuite) Saint Michel de Bruxelles vers 1910 \cite{Lambert}, pages 43-44. Un autre élève de Mansion joue un rôle important dans la formation mathématique de Lemaître, Ernest Pasquier (1849-1926), professeur à l'Université Catholique de Louvain.

Un des élèves les plus brillants de Paul Mansion est  Alphonse Demoulin\footnote{ En 1927, Alphonse Demoulin prononce l'éloge de Paul Mansion devant la Classe des Sciences de l'Académie royale des sciences, des lettres et des beaux-arts de Belgique\cite{Demoulin}.}. Demoulin a  tout d'abord Paul Mansion à l'Ecole Normale des sciences. En Octobre 1893, après son retour d'un séjour effectué à Leipzig\footnote{Voir le témoignage de Demoulin sur les enseignements qu'il a suivis à Paris puis Leipzig \cite{Lef3}.} dans l'Institut de Mathématiques animé par Sophus Lie, Alphonse Demoulin commence à enseigner à l'Université de Gand. Il y restera toute sa carrière mais cessera ses enseignements en 1936 lors de la flamandisation de l'Université de Gand. Les travaux mathématiques d'Alphonse Demoulin se situent essentiellement dans le domaine de la Géométrie Différentielle\cite{Maw1}. Les recherches de Demoulin sont récompensées à plusieurs reprises par l'Académie des Sciences de Paris : il reçoit le prix du baron de Joest en 1906, le prix Bordin en 1911\footnote{Le sujet du concours proposé était : \og {\it perfectionner en un point important la théorie des systèmes triples de surfaces orthogonales} \fg.} et le prix Poncelet en 1945. En Belgique, il se voit décerner en 1919 Ie prix décennal des mathématiques pures pour la période 1904-1913.

En 1882, Paul Mansion est élu membre correspondant de la Classe des Sciences de l'Académie Royale de Belgique et il en devient membre titulaire en 1887. Paul Mansion adresse une lettre\footnote{Dossier Mansion de l'Académie Royale de Belgique.} le 18 Décembre 1887 au général Liagre, secrétaire perpétuel de l'Académie :
\begin{quote}
\begin{it}
Veuillez, je vous prie, être mon interprète auprès de la Classe des Sciences de l'Académie royale de Belgique pour lui transmettre mes bien vifs remerciements à l'occasion de mon élection comme membre titulaire, dans sa dernière séance.
\end{it}
\end{quote}
Paul Mansion préside durant l'année 1889-1890 la Société Scientifique de Bruxelles dont il est membre depuis sa création en 1875. Un des créateurs de cette société savante est le père jésuite, et scientifique, belge Ignace Carbonnelle (1829-1889). En 1875, Ignace Carbonnelle enseigne au Collège Saint Michel de Bruxelles. Plus tard, en 1907, Paul Mansion, alors professeur émérite de l'Université de Gand est le secrétaire de la Société Scientifique de Bruxelles. Cette société est organisée en {\it sections}. La première section est celle des {\it Math\'ematiques, Astronomie, Géodésie et Mécanique}, présidée par Georges Humbert (1859-1921). Charles de La Vallée Poussin (1866-1962) et Robert d'Adhémar (1874-1941) en assurent la vice-présidence et Hector Dutordoir\footnote{Hector Dutordoir est un ingénieur en chef, directeur du service technique provincial, demeurant à Gand.} le secrétariat. La Société Scientifique de Bruxelles est très active et publie régulièrement ses actes dans {\it Annales de la Société Scientifique de Bruxelles}.

\section{Sur les écrits de Paul Mansion en probabilités}
 
Paul Mansion a écrit les articles suivants en Probabilités :
\begin{itemize}
\item Sur le problème des partis. Mém. in-8, 1868. T. XXI, 13 pp, mémoire Acad. Sci. Bel.
\item Note sur la portée objective du calcul des probabilités. En appendice, dans
Paul Janet : Les causes finales, pp. 720-725. Paris, G. Baillière, 1882; sous le titre :
L'argument épicurien et le calcul des probabilités.
\item Sur le théorème de Jacques Bernoulli. 1892. T. XVI, 1e partie, pp. 85-87, Ann. Soc. Sci. Bruxelles.
\item Sur la loi des grands nombres de Poisson. 1893. T. XXV, pp. 11-13, Bulletin (3e série).
\item Sur la loi des grands nombres de Poisson. 1892-1893. T. XXII, p. 56. Messenger of Mathematics ; New series.
\item Sur le théorème de Jacques Bernoulli. 1898. T. XXII, 1e partie, pp. 3-4.
\item Démonstration du théorème de Jacques Bernoulli. Congrès de Munich, 1900;
Fribourg, Herder, 1901 ; pp. 427-428, Congrès scientifique international des catholiques.
\item Démonstration du théorème de Bernoulli. Sur une intégrale considérée en
calcul des probabilités. Sur une intégrale considérée par Poisson en calcul des probabilités.
1902. T. XXVI, 2e partie, pp. 191-214; 2e partie, p. 126, Ann. Soc. Sci. Bruxelles.
\item Sur la portée objective du calcul des probabilités. 1903, pp. 1235-1294, Bulletin de la Classe des sciences.
\item Sur la loi des grands nombres de Poisson. Sur une sommation d'intégrales
considérées en calcul des probabilités. 1904. T. XXVIII, 1e partie, pp. 72-77, 166-167, Ann. Soc. Sci. Bruxelles.
\item Démonstration de la loi des grands nombres. 1910, pp. 158-160. Bulletin de la Classe des sciences.
\item Démonstration de la loi des grands nombres. 1910. T. XXXIV, pp. 230-233, Ann. Soc. Sci. Bruxelles.
\end{itemize}

Quel est le contenu du premier article de Mansion ?  En reprenant l'analyse et les notations de \cite{PG} , Paul Mansion dans son mémoire présenté le 12 Mai 1868 à l'Académie Royale des Sciences de Belgique, considère une distribution \og binomiale négative\fg\ (ou {\it binomial waiting time distribution}) :
\begin{equation}
P_{s_{1},s_{2}}^{(1)}=p_{1}^{s_{1}}\sum_{j=0}^{s_{2}-1}
\left(
\begin{array}{c}
s_{1}-1+j\\
s_{1}-1
\end{array}\right)
(1-p_{1})^{j}.
\end{equation}
Mansion établit alors pour $2$ joueurs ( et plus) :
\begin{equation}
P_{s_{1},s_{2}}^{(1)}=\frac{\Gamma(s_{1}+s_{2})}{\Gamma(s_{1})\Gamma(s_{2})}
\sum_{i=0}^{s_{2}-1}(-1)^{i}\frac{1}{s_{1}+i}
\left(
\begin{array}{c}
s_{2}-1\\
i
\end{array}\right)
p_{1}^{i+s_{1}}.
\end{equation}
En 1878 Catalan propose dans \cite{Catalan1} une solution et différentes interprétations au problème des partis en partant lui aussi d'une distribution binomiale négative (voir \cite{PG}).

En 1895, Edouard Goedseels (1857-1928)\footnote{Edouard Goedseels est membre fondateur de la soci\'et\'e Scientifique de Bruxelles et administrateur-inspecteur de l'Observatoire royal de Belgique. Il \'etait aussi professeur \`a l'Universit\'e Catholique de Louvain \cite{Maw2}} publie dans les {\it Annales de la société scientifique de Bruxelles}, tome XVII p. 8,  une démonstration  théorème de Bernoulli. Les deux articles de Paul Mansion sur le th\'eor\`eme de Bernoulli paraissent en 1902 et 1904 dans la continuité donc des travaux de Goedseels. Plus tard, en 1907, Charles de La Vallée Poussin publie dans  la continuité des travaux de Goedseels et de Mansion, deux articles sur le théorème de Bernoulli. Les trois math\'ematiciens cherchent à préciser la vitesse de convergence de la probabilit\'e $P\left(\left\vert\frac{S_{n}}{n}-p\right\vert\leq \epsilon\right)$ vers $1$, $n\to\infty$, $S_{n}$ \'etant une somme de variables de Bernoulli ind\'ependantes de probabilité $p$ de succès. Paul Mansion souligne la difficult\'e de l'utilisation du th\'eor\`eme de Bernoulli, dans sa forme originale, dans les applications pratiques. Les travaux de Mansion et de La Vallée Poussin, pour reprendre une expression de Seneta, \og présagent d'un retour aux méthodes exactes\fg. Le mathématicien français Robert de Montessus de Ballore, professeur à l'Université Catholique de Lille, cite ces différents auteurs dans son livre sur les probabilités publié en 1908 \cite{RMB1908}. On consultera aussi \cite{Lef1} pour plus d'éléments.

En lien avec le calcul des probabilités, Paul Mansion s'est toujours  intéressé au développement du Calcul Intégral. Ainsi,par exemple, Il publie en 1885 dans {\it Mathesis} une preuve de la seconde formule de moyenne\footnote{Ce résultat est surtout utilisé dans la théorie des séries de Fourier.} que Léopold Kronecker (1823-1891) lui a transmise (voir \cite{Lef2} pour une analyse de l'article). Dans le fonds Mansion de la Bibliothèque Royale se trouvent justement les lettres de Kronecker contenant la démonstration.

\section{Conclusion}
L'oeuvre scientifique de Paul Mansion est immense. Par la qualité de ses travaux mathématiques, mais aussi par ses activités éditoriales et par responsabilités à l'Académie Royale et à la Société Scientifique de Bruxelles, Paul Mansion a eu une grande influence sur les mathématiques belges de son temps, en Recherche mais aussi dans l'Enseignement. Paul Mansion était en contact avec de nombreux mathématiciens et non des moindre d'ailleurs. 

Nous avons évoqué les travaux en probabilités de Paul Mansion héritier d'une tradition belge dans ce domaine.  Ces travaux ne constituent qu'une très petite part de ses écrits. Les interlocuteurs de Paul Mansion sur ce sujet furent essentiellement dans les années 1900 les mathématiciens belges Goedseels et De La Vallée Poussin\footnote{De La Vallée a lui aussi très peu écrit dans le domaine des Probabilités.} et le mathématicien français Robert de Montessus. On pourra ajouter le mathématicen et astronome de l'Université de Louvain, Maurice Alliaume (1882-1931)\footnote{Il faudrait d'alleurs revenir à la place que tient  l'Astronomie comme source de sujets pour les probabilités et statistiques. Le travaux de Goedseels (et sa position à l'Observatoire Royal) sur la méthodes des moindres carrés sont significatifs.}. Robert de Montessus dans son livre de leçons de probabilités paru en 1908 cite abondamment les travaux de ces mathématiciens belges. Au sujet de ce livre, une question se pose. Robert de Montessus dédie un chapitre de son ouvrage aux travaux novateurs de Louis Bachelier (1870-1946) sur la théorie de la spéculation. Les mathématiciens belges en ce début du XXe siècle, Mansion en premier lieu, se sont-ils intéressés aux idées de Bachelier ?

Terminons par une remarque sur le théorème des nombres premiers démontré de façon indépendante en 1896 par Charles de La Vallée Poussin et Jacques Hadamard (1865-1963). Dans \cite{Maw3}, Jean Jean Mawhin indique que comme Hadamard, De La Vallée Poussin n'évoque pas Karl Friedrich Gauss (1777-1855) dans son article contenant sa démonstration du théorème des nombres premiers . Mais en 1899, dans un article sur la fonction $\zeta$ de Riemann,  De La Vallée Poussin cite cette fois-ci Gauss . C'est Paul Mansion qui a rédigé le rapport de cet article pour publication dans les annales de l'Académie Royale.

%Peut-on à partir de ces quelques éléments voir apparaître un réseau de mathématiciens dans le %domaine des probabilités avec Paul Mansion en son centre ? Si oui, doit-on se limiter à un réseau %uniquement belge, hormis Robert de Montessus qui est à Lille ?
\newpage

\begin{tikzpicture}
\draw (4.4,0) node[below]{$1844$};
\draw (11.9,0) node[below]{$1919$};
\draw [line width=2pt](4.4,0) --(11.9,0) node[right]{Mansion};
%\draw [line width=5pt] (-3.4,0)--(4,0) node[above left=5pt]{Garnier};
\draw[line width=2pt] (-3.4,0)--(4,0) node[above left=5pt]{Garnier};
\draw (-3.4,0) node[above]{$1766$};
\draw (4,0) node[above]{$1840$};
%\draw (0,0) circle (1) ;
\end{tikzpicture}

\begin{tikzpicture}
\draw (4.4,0) node[below]{$1844$};
\draw (11.9,0) node[below]{$1919$};
\draw[line width=2pt] (4.4,0) --(11.9,0) node[right]{Mansion};
\draw[line width=2pt] (-0.5,0)--(8.4,0) node[above left=5pt]{Manderlier};
\draw[line width=2pt] [color=red] (4.4,0)--(8.4,0) ;
\draw (-0.5,0) node[above]{$1795$};
\draw (8.4,0) node[above]{$1884$};
%\draw (0,0) circle (1) ;
\end{tikzpicture}

\begin{tikzpicture}
\draw (4.4,0) node[below]{$1844$};
\draw (11.9,0) node[below]{$1919$};
\draw [line width=2pt](4.4,0) --(11.9,0) node[right]{Mansion};
\draw [line width=2pt] (-0.4,0)--(7.4,0) node[above left=5pt]{Quetelet};
\draw[line width=2pt] [color=red] (4.4,0)--(7.4,0) ;
\draw (-0.4,0) node[above]{$1796$};
\draw (7.4,0) node[above]{$1874$};
%\draw (0,0) circle (1) ;
\end{tikzpicture}

\begin{tikzpicture}
\draw (4.4,0) node[below]{$1844$};
\draw (11.9,0) node[below]{$1919$};
\draw [line width=2pt](4.4,0) --(11.9,0) node[right]{Mansion};
\draw [line width=2pt] (-0.3,0)--(5.2,0) node[above left=5pt]{Lemaire};
\draw[line width=2pt] [color=red] (4.4,0)--(5.2,0) ;
\draw (-0.3,0) node[above]{$1797$};
\draw (5.2,0) node[above]{$1852$};
%\draw (0,0) circle (1) ;
\end{tikzpicture}

\begin{tikzpicture}
\draw (4.4,0) node[below]{$1844$};
\draw (11.9,0) node[below]{$1919$};
\draw [line width=2pt](4.4,0) --(11.9,0) node[right]{Mansion};
\draw [line width=2pt] (0.1,0)--(6.4,0) node[above left=5pt]{Timmermans};
\draw[line width=2pt] [color=red] (4.4,0)--(6.4,0) ;
\draw (0.1,0) node[above]{$1801$};
\draw (6.4,0) node[above]{$1864$};
%\draw (0,0) circle (1) ;
\end{tikzpicture}

\begin{tikzpicture}
\draw (4.4,0) node[below]{$1844$};
\draw (11.9,0) node[below]{$1919$};
\draw [line width=2pt](4.4,0) --(11.9,0) node[right]{Mansion};
\draw [line width=2pt] (0.3,0)--(5.7,0) node[above left=5pt]{Meyer};
\draw[line width=2pt] [color=red] (4.4,0)--(5.7,0) ;
\draw (0.3,0) node[above]{$1803$};
\draw (5.7,0) node[above]{$1857$};
%\draw (0,0) circle (1) ;
\end{tikzpicture}

\begin{tikzpicture}
\draw (4.4,0) node[below]{$1844$};
\draw (11.9,0) node[below]{$1919$};
\draw[line width=2pt] (4.4,0) --(11.9,0) node[right]{Mansion};
\draw [line width=2pt] (0.4,0)--(4.9,0) node[above left=5pt]{Verhulst};
\draw[line width=2pt] [color=red] (4.4,0)--(4.9,0) ;
\draw (0.4,0) node[above]{$1804$};
\draw (4.9,0) node[above]{$1849$};
%\draw (0,0) circle (1) ;
\end{tikzpicture}

\begin{tikzpicture}
\draw (4.4,0) node[below]{$1844$};
\draw (11.9,0) node[below]{$1919$};
\draw [line width=2pt](4.4,0) --(11.9,0) node[right]{Mansion};
\draw [line width=2pt] (1.4,0)--(9.4,0) node[above left=5pt]{Catalan};
\draw[line width=2pt] [color=red] (4.4,0)--(9.4,0) ;
\draw (1.4,0) node[above]{$1814$};
\draw (9.4,0) node[above]{$1894$};
%\draw (0,0) circle (1) ;
\end{tikzpicture}

\begin{tikzpicture}
\draw (4.4,0) node[below]{$1844$};
\draw (11.9,0) node[below]{$1919$};
\draw [line width=2pt](4.4,0) --(11.9,0) node[right]{Mansion};
\draw [line width=2pt] (1.7,0)--(6.7,0) node[above left=5pt]{Schaar};
\draw[line width=2pt] [color=red] (4.4,0)--(6.7,0) ;
\draw (1.7,0) node[above]{$1817$};
\draw (6.7,0) node[above]{$1867$};
%\draw (0,0) circle (1) ;
\end{tikzpicture}

\begin{tikzpicture}
\draw (4.4,0) node[below]{$1844$};
\draw (11.9,0) node[below]{$1919$};
\draw [line width=2pt](4.4,0) --(11.9,0) node[right]{Mansion};
\draw [line width=2pt] (2,0)--(9.3,0) node[above left=5pt]{Boudin};
\draw[line width=2pt] [color=red] (4.4,0)--(9.3,0) ;
\draw (2,0) node[above]{$1820$};
\draw (9.3,0) node[above]{$1893$};
%\draw (0,0) circle (1) ;
\end{tikzpicture}

\begin{tikzpicture}
\draw (4.4,0) node[below]{$1844$};
\draw (11.9,0) node[below]{$1919$};
\draw [line width=2pt](4.4,0) --(11.9,0) node[right]{Mansion};
\draw [line width=2pt] (2.9,0)--(8.9,0) node[above left=5pt]{Dauge};
\draw[line width=2pt] [color=red] (4.4,0)--(8.9,0) ;
\draw (2.9,0) node[above]{$1829$};
\draw (8.9,0) node[above]{$1889$};
%\draw (0,0) circle (1) ;
\end{tikzpicture}

\begin{tikzpicture}
\draw (4.4,0) node[below]{$1844$};
\draw (11.9,0) node[below]{$1919$};
\draw [line width=2pt]node[right]{Mansion}(4.4,0) --(11.9,0) ;
\draw [line width=2pt] (4.9,0)--(12.6,0) node[above left=5pt]{Pasquier};
\draw[line width=2pt] [color=red] (4.9,0)--(11.9,0) ;
\draw (4.9,0) node[above]{$1849$};
\draw (12.6,0) node[above]{$1926$};
%\draw (0,0) circle (1) ;
\end{tikzpicture}

\begin{tikzpicture}
\draw (4.4,0) node[below]{$1844$};
\draw (11.9,0) node[below]{$1919$};
\draw [line width=2pt]node[right]{Mansion}(4.4,0) --(11.9,0) ;
\draw [line width=2pt] (5.2,0)--(12.8,0) node[above left=5pt]{Bosmans};
\draw[line width=2pt] [color=red] (5.2,0)--(11.9,0) ;
\draw (5.2,0) node[above]{$1852$};
\draw (12.8,0) node[above]{$1928$};
%\draw (0,0) circle (1) ;
\end{tikzpicture}

\begin{tikzpicture}
\draw (4.4,0) node[below]{$1844$};
\draw (11.9,0) node[below]{$1919$};
\draw [line width=2pt]node[right]{Mansion} (4.4,0) --(11.9,0) ;
\draw [line width=2pt] (7.5,0)--(14.9,0) node[above left=5pt]{Claeys};
\draw[line width=2pt] [color=red] (7.5,0)--(11.9,0) ;
\draw (7.5,0) node[above]{$1875$};
\draw (14.9,0) node[above]{$1949$};
%\draw (0,0) circle (1) ;
\end{tikzpicture}

\begin{tikzpicture}
\draw (4.4,0) node[below]{$1844$};
\draw (11.9,0) node[below]{$1919$};
\draw [line width=2pt]node[right]{Mansion} (4.4,0) --(11.9,0) ;
\draw [line width=2pt] (8.2,0)--(13.1,0) node[above left=5pt]{Alliaume};
\draw[line width=2pt] [color=red] (8.2,0)--(11.9,0) ;
\draw (8.2,0) node[above]{$1882$};
\draw (13.1,0) node[above]{$1931$};
%\draw (0,0) circle (1) ;
\end{tikzpicture}

\begin{tikzpicture}
\draw (4.4,0) node[below]{$1844$};
\draw (11.9,0) node[below]{$1919$};
\draw [line width=2pt] node[right]{Mansion}(4.4,0) --(11.9,0) ;
\draw [line width=2pt] (8.4,0)--(15.6,0) node[above left=5pt]{Sarton};
\draw[line width=2pt] [color=red] (8.4,0)--(11.9,0) ;
\draw (8.4,0) node[above]{$1884$};
\draw (15.6,0) node[above]{$1956$};
%\draw (0,0) circle (1) ;
\end{tikzpicture}

\end{document}